\documentclass[10pt]{amsart}
\usepackage{amsfonts}
\usepackage{amsmath}
\usepackage{amsthm}
\usepackage{amssymb}
\usepackage{latexsym}
\newtheorem{cor}{Corollary}[section]
\newtheorem{lem}{Lemma}[section]

\newtheorem{prop}{Proposition}[section]

\theoremstyle{definition}
\newtheorem{defn}{Definition}[section]
\theoremstyle{definition}
\newtheorem{thm}{Theorem}

\newtheorem*{rem}{Remark}

\newenvironment{pf}{\proof}{\endproof}

\theoremstyle{remark}

\numberwithin{equation}{section}
\begin{document}

\newcommand{\thmref}[1]{Theorem~\ref{#1}}
\newcommand{\secref}[1]{Sect.~\ref{#1}}
\newcommand{\lemref}[1]{Lemma~\ref{#1}}
\newcommand{\propref}[1]{Proposition~\ref{#1}}
\newcommand{\corref}[1]{Corollary~\ref{#1}}
\newcommand{\remref}[1]{Remark~\ref{#1}}
\newcommand{\nc}{\newcommand}
\newcommand{\rnc}{\renewcommand}
\nc{\cal}{\mathcal}
\nc{\goth}{\mathfrak}
\rnc{\bold}{\mathbf}
\renewcommand{\frak}{\mathfrak}
\renewcommand{\Bbb}{\mathbb}

\nc{\Cal}{\mathcal}
\nc{\Xp}[1]{X^+(#1)}
\nc{\Xm}[1]{X^-(#1)}
\nc{\on}{\operatorname}
\nc{\ch}{\mbox{ch}}
\nc{\Z}{{\bold Z}}
\nc{\J}{{\cal J}}
\nc{\C}{{\bold C}}
\nc{\Q}{{\bold Q}}
\renewcommand{\P}{{\cal P}}
\nc{\N}{{\Bbb N}}
\nc\beq{\begin{equation}}
\nc\enq{\end{equation}}
\nc\lan{\langle}
\nc\ran{\rangle}
\nc\bsl{\backslash}
\nc\mto{\mapsto}
\nc\lra{\leftrightarrow}
\nc\hra{\hookrightarrow}
\nc\sm{\smallmatrix}
\nc\esm{\endsmallmatrix}
\nc\sub{\subset}
\nc\ti{\tilde}
\nc\nl{\newline}
\nc\fra{\frac}
\nc\und{\underline}
\nc\ov{\overline}
\nc\ot{\otimes}
\nc\bbq{\bar{\bq}_l}
\nc\bcc{\thickfracwithdelims[]\thickness0}
\nc\ad{\text{\rm ad}}
\nc\Ad{\text{\rm Ad}}
\nc\Hom{\text{\rm Hom}}
\nc\End{\text{\rm End}}
\nc\Ind{\text{\rm Ind}}
\nc\Res{\text{\rm Res}}
\nc\Ker{\text{\rm Ker}}
\rnc\Im{\text{Im}}
\nc\sgn{\text{\rm sgn}}
\nc\tr{\text{\rm tr}}
\nc\Tr{\text{\rm Tr}}
\nc\supp{\text{\rm supp}}
\nc\card{\text{\rm card}}
\nc\bst{{}^\bigstar\!}
\nc\he{\heartsuit}
\nc\clu{\clubsuit}
\nc\spa{\spadesuit}
\nc\di{\diamond}

\nc\al{\alpha}
\nc\bet{\beta}
\nc\ga{\gamma}
\nc\de{\delta}
\nc\ep{\epsilon}
\nc\io{\iota}
\nc\om{\omega}
\nc\si{\sigma}
\rnc\th{\theta}
\nc\ka{\kappa}
\nc\la{\lambda}
\nc\ze{\zeta}

\nc\vp{\varpi}
\nc\vt{\vartheta}
\nc\vr{\varrho}

\nc\Ga{\Gamma}
\nc\De{\Delta}
\nc\Om{\Omega}
\nc\Si{\Sigma}
\nc\Th{\Theta}
\nc\La{\Lambda}
\nc\boa{\bold a}
\nc\bob{\bold b}
\nc\boc{\bold c}
\nc\bod{\bold d}
\nc\boe{\bold e}
\nc\bof{\bold f}
\nc\bog{\bold g}
\nc\boh{\bold h}
\nc\boi{\bold i}
\nc\boj{\bold j}
\nc\bok{\bold k}
\nc\bol{\bold l}
\nc\bom{\bold m}
\nc\bon{\bold n}
\nc\boo{\bold o}
\nc\bop{\bold p}
\nc\boq{\bold q}
\nc\bor{\bold r}
\nc\bos{\bold s}
\nc\bou{\bold u}
\nc\bov{\bold v}
\nc\bow{\bold w}
\nc\boz{\bold z}

\nc\ba{\bold A}
\nc\bb{\bold B}
\nc\bc{\bold C}
\nc\bd{\bold D}
\nc\be{\bold E}
\nc\bg{\bold G}
\nc\bh{\bold h}
\nc\bH{\bold H}

\nc\bi{\bold I}
\nc\bj{\bold J}
\nc\bk{\bold K}
\nc\bl{\bold L}
\nc\bm{\bold M}
\nc\bn{\bold N}
\nc\bo{\bold O}
\nc\bp{\bold P}
\nc\bq{\bold Q}
\nc\br{\bold R}
\nc\bs{\bold S}
\nc\bt{\bold T}
\nc\bu{\bold U}
\nc\bv{\bold v}
\nc\bV{\bold V}

\nc\bw{\bold W}
\nc\bz{\bold Z}
\nc\bx{\bold x}
\nc\bX{\bold X}
\nc\blambda{{\mbox{\boldmath $\Lambda$}}}
\nc\bpi{{\mbox{\boldmath $\pi$}}}

\nc\e[1]{E_{#1}}
\nc\ei[1]{E_{\delta - \alpha_{#1}}}
\nc\esi[1]{E_{s \delta - \alpha_{#1}}}
\nc\eri[1]{E_{r \delta - \alpha_{#1}}}
\nc\ed[2][]{E_{#1 \delta,#2}}
\nc\ekd[1]{E_{k \delta,#1}}
\nc\emd[1]{E_{m \delta,#1}}
\nc\erd[1]{E_{r \delta,#1}}

\nc\ef[1]{F_{#1}}
\nc\efi[1]{F_{\delta - \alpha_{#1}}}
\nc\efsi[1]{F_{s \delta - \alpha_{#1}}}
\nc\efri[1]{F_{r \delta - \alpha_{#1}}}
\nc\efd[2][]{F_{#1 \delta,#2}}
\nc\efkd[1]{F_{k \delta,#1}}
\nc\efmd[1]{F_{m \delta,#1}}
\nc\efrd[1]{F_{r \delta,#1}}
\nc{\ug}{\bu^{fin}}

\nc\fa{\frak a}
\nc\fb{\frak b}
\nc\fc{\frak c}
\nc\fd{\frak d}
\nc\fe{\frak e}
\nc\ff{\frak f}
\nc\fg{\frak g}
\nc\fh{\frak h}
\nc\fj{\frak j}
\nc\fk{\frak k}
\nc\fl{\frak l}
\nc\fm{\frak m}
\nc\fn{\frak n}
\nc\fo{\frak o}
\nc\fp{\frak p}
\nc\fq{\frak q}
\nc\fr{\frak r}
\nc\fs{\frak s}
\nc\ft{\frak t}
\nc\fu{\frak u}
\nc\fv{\frak v}
\nc\fz{\frak z}
\nc\fx{\frak x}
\nc\fy{\frak y}

\nc\fA{\frak A}
\nc\fB{\frak B}
\nc\fC{\frak C}
\nc\fD{\frak D}
\nc\fE{\frak E}
\nc\fF{\frak F}
\nc\fG{\frak G}
\nc\fH{\frak H}
\nc\fJ{\frak J}
\nc\fK{\frak K}
\nc\fL{\frak L}
\nc\fM{\frak M}
\nc\fN{\frak N}
\nc\fO{\frak O}
\nc\fP{\frak P}
\nc\fQ{\frak Q}
\nc\fR{\frak R}
\nc\fS{\frak S}
\nc\fT{\frak T}
\nc\fU{\frak U}
\nc\fV{\frak V}
\nc\fZ{\frak Z}
\nc\fX{\frak X}
\nc\fY{\frak Y}
\nc\tfi{\ti{\Phi}}
\nc\bF{\bold F}

\nc\ua{\bold U_\A}

%%%%%%%%%%%%%%%%%%%%%%%%%%%%%%%%%%%%%%%%%%%%%%%%%%%%%%
\nc\qinti[1]{[#1]_i}
\nc\q[1]{[#1]_q}
\nc\xpm[2]{E_{#2 \delta \pm \alpha_#1}}  %\xpm{j}{l}
\nc\xmp[2]{E_{#2 \delta \mp \alpha_#1}}
\nc\xp[2]{E_{#2 \delta + \alpha_{#1}}}
\nc\xm[2]{E_{#2 \delta - \alpha_{#1}}}
\nc\hik{\ed{k}{i}}
\nc\hjl{\ed{l}{j}}
\nc\qcoeff[3]{\left[ \begin{smallmatrix} {#1}& \\ {#2}& \end{smallmatrix}
\negthickspace \right]_{#3}}
\nc\qi{q}
\nc\qj{q}

\nc\ufdm{{_\ca\bu}_{\rm fd}^{\le 0}}

%%%%%%%%%%%%%%%%%%%%%%%%%%%%%%%%%%%%%%%%%%%%%%%%%%%%%%

%\nc\rtimes
\nc\isom{\cong}

\nc{\pone}{{\Bbb C}{\Bbb P}^1}
\nc{\pa}{\partial}
\def\H{\cal H}
\def\L{\cal L}
\nc{\F}{{\cal F}}
\nc{\Sym}{{\goth S}}
\nc{\A}{{\cal A}}
\nc{\arr}{\rightarrow}
\nc{\larr}{\longrightarrow}

\nc{\ri}{\rangle}
\nc{\lef}{\langle}
\nc{\W}{{\cal W}}
\nc{\uqatwoatone}{{U_{q,1}}(\su)}
\nc{\uqtwo}{U_q(\goth{sl}_2)}
\nc{\dij}{\delta_{ij}}
\nc{\divei}{E_{\alpha_i}^{(n)}}
\nc{\divfi}{F_{\alpha_i}^{(n)}}
\nc{\Lzero}{\Lambda_0}
\nc{\Lone}{\Lambda_1}
\nc{\ve}{\varepsilon}
\nc{\phioneminusi}{\Phi^{(1-i,i)}}
\nc{\phioneminusistar}{\Phi^{* (1-i,i)}}
\nc{\phii}{\Phi^{(i,1-i)}}
\nc{\Li}{\Lambda_i}
\nc{\Loneminusi}{\Lambda_{1-i}}
\nc{\vtimesz}{v_\ve \otimes z^m}

\nc{\asltwo}{\widehat{\goth{sl}_2}}
\nc\eh{\frak h^e}
\nc\loopg{L(\frak g)}
\nc\eloopg{L^e(\frak g)}
\nc\ebu{\bu^e}
\nc\loopa{L(\frak a)}

\nc\teb{\tilde E_\boc}
\nc\tebp{\tilde E_{\boc'}}

\title{Braid group actions and tensor products}\author{Vyjayanthi Chari}
\address{Vyjayanthi Chari, Department of Mathematics, University of California,
Riverside, CA 92521.}
\maketitle

{\small{\centerline {Abstract}

{\it \noindent We define an action of the braid group of a simple
Lie algebra on the space of  imaginary roots in the corresponding
quantum affine algebra. We then use this action to determine an
explicit  condition for  a tensor product of arbitrary irreducible
finite--dimensional representations is cyclic. This allows us to
determine the  set of points at which the corresponding
$R$--matrix has a zero.}}

\section{Introduction} In this paper we  give a sufficient condition for the
tensor product of irreducible finite--dimensional representations
of quantum affine algebras to be cyclic. This condition is
obtained by defining a braid group action on the imaginary root
vectors. We  make the condition  explicit in  Section 5 and see
that it is a natural generalization of the condition in
\cite{CPqa} given in the $sl_2$ case. This allows us for instance,
to determine the finite set of points at which a tensor product of
fundamental representations can fail to be cyclic.
 Our result proves a generalization of a recent result of
Kashiwara \cite{K}, \cite{VV}, \cite{Nak1}. Further, it also
establishes a conjecture stated in \cite{K}, \cite{HKOTY}.

We describe our results in some detail. Let ${\frak{g}}$ be a    complex simple
finite--dimensional Lie algebra of rank $n$, and let $\bu_q$ be the quantized
untwisted affine algebra over $\bc(q)$ associated to $\frak g$. For every
$n$--tuple $\bpi=(\pi_1,\cdots ,\pi_n)$ of polynomials  with coefficients in
$\bc(q)[u]$ and with constant term one, there exists  a unique (up to
isomorphism)  irreducible finite--dimensional representation  $V(\bpi)$ of
$\bu_q$. For each element  $w$ in the Weyl group $W$ of $\frak g$,  let
$v_{w\bpi}$ be the extremal vector defined in \cite{K}.
 In this paper we compute the action of the imaginary root vectors in $\bu_q$ on
the elements $v_{w\bpi}$. To do this we define in Section 2 an
action of the braid group $\cal{B}$  of $\frak g$ on elements of
$(\bc(q)[[u]])^n$ and prove that the eigenvalue of $v_{w\bpi}$ is
the element $T_w(\bpi)$ where $T:W\to\cal{B}$ is the canonical
section defined in \cite{Bo}.

Let $\pi(u)\in\bc(q)[u]$ be a polynomial that splits in $\bc(q)$. Any such polynomial $\pi$
can be written uniquely as a product
\begin{equation*}\pi(u)=\prod_{r=1}^k(1-a_rq^{m_r-1}u)(1-a_rq^{m_r-3}u)\cdots
(1-a_rq^{-m_r+1}u),\end{equation*}
where $a_r\in\bc(q)$ and $m_r\in\bz_+$ satisfy
\begin{equation*}  \frac{a_r}{a_l} \ne  q^{\pm(m_r+ m_l-2m)},\ \ 0\le
m<\text{min}(m_r,m_l), \ \ \end{equation*}  if $r<l$. Let $S(\pi)$
be the collection of the pairs $(a_r,m_r)$. $1\le r\le k$ defined
above. Say that a polynomial $\pi'(u)$ is in general  position
with respect to $\pi(u)$, if
\begin{equation*}\frac{a_r'}{a_\ell}\ne  q^{-(m_{r}'+m_\ell-2p)},  0\le p<  m_{r},\end{equation*}
for all pairs $(a_r',m_r')\in S(\pi')$ and $(a_\ell, m_\ell)\in S(\pi)$.

\medspace Our main result Theorem \ref{main} is the following, we
restrict ourselves to the simply laced (only in the introduction).
Let $s_1, s_2, \cdots ,s_n$ be the set of simple reflections in
$W$ and let $s_{i_1}\cdots s_{i_N}$ be  a reduced expression for
the longest element $w_0\in W$.

 \noindent {\it{ The tensor product
$V(\bpi')\otimes V(\bpi)$ is cyclic on $v_{\bpi'}\otimes v_{\bpi}$
if for all $1\le j\le N$ the polynomial $(T_{i_{j+1}}\cdots
T_{i_N}\bpi')_{i_{j}}$ is in general position with  respect to
$\bpi_{i_j}$.
 More generally, let
$V_1,\cdots V_r$ be  irreducible finite--dimensional
representations such that   $V_j\otimes V_l$ is cyclic if $j<\ell
$. Then $V_1\otimes\cdots \otimes V_r$ is cyclic. }}

\noindent Of course, we first have to prove that
$(T_{i_{j+1}}\cdots T_{i_N}\bpi')_{i_{j}}$ is a polynomial, we do
this in Section 2.  In Section 5, we write down the polynomials
$(T_{i_{j+1}}\cdots T_{i_N}\bpi')_{i_{j}}$, for all $\frak g$ and
a specific reduced expression of $w_0$ thus making the condition
for a tensor product to be cyclic explicit and we see that it is
the appropriate generalization of the result in the case when
$\frak g =sl_2$.

To make the connection with Kashiwara's conjectures, we consider
the
 case
\begin{equation*} \pi_j(u) =1\ \ (j\ne i), \ \
\pi_i(u)=\prod_{s=1}^m(1-q^{m+1-2s}au)\ \ (a\in\bc(q)).
\end{equation*} Denoting this $n$--tuple of polynomials as
$\bpi^i_{m,a}$ and the corresponding representation by
$V(\bpi^i_{m,a}$, we prove the following: \medspace \noindent
{\it{ Let $l\ge 1$ and let $i_j\in I$, $m_j\in\bz_+$,
$a_j\in\bc(q)$ for $1\le j\le l$. The tensor product
$V(\bpi^{i_1}_{m_1,a_1})\otimes V(\bpi^{i_2}_{m_2,
a_2})\otimes\cdots\otimes V(\bpi^{i_l}_{m_l,a_l})$ is cyclic on
the tensor product of highest weight vectors if for all $r<s$,
\begin{equation*} \frac{a_r}{a_s}\ne q^{m_r-m_s-p}, \ \ \forall \ p\ge
0.\end{equation*}}} In  Corollary 5.1, we write down the precise
values of $a_r/a_s$ at which the tensor product is not cyclic. In
the case of two when $\ell =2$, this
 is  the set of all possible zeros of the $R$--matrix, $R(a):
V(\bpi^{i_1}_{m_1,a})\otimes V(\bpi^{i_2}_{m_2, 1})\to
V(\bpi^{i_1}_{m_1,a})\otimes V(\bpi^{i_2}_{m_2, 1})$.

\noindent The case when $m_j=1$ was originally conjectured  and
partially proved in \cite{AK} and  completely proved in \cite{K}
(and in \cite{VV} for the simply-laced case),
\cite{Nak1},\cite{Nak2}. The result in the case  when the $m_i$
are arbitrary but $a_i=1$  for all $i$ was conjectured in
\cite{K}, \cite{HKOTY}. In the case of $A_n$ and $C_n$ and when
$m_j=1$ the values of $a_r/a_s$ where the tensor product is not
cyclic was written in \cite{AK}

 In the exceptional case
 the explicit calculation of the possible points of
reducibility  can be used  to write down the precise $\frak
g$--module structures of the exceptional algebras for {\it all}
nodes of the Dynkin diagram. Details of this will appear
elsewhere.

Finally, recall that a tensor product of two irreducible
finite--dimensional representations is irreducible if both
$V\otimes V'$ and its dual are cyclic on the tensor product of
highest weight vectors. Thus our theorem gives us a sufficient
condition for the tensor product $V\otimes V'$ to be irreducible.
When $\frak g=sl_2$, this condition is the same as the one given
in \cite{CPqa}.

\section{Preliminaries} In this section we recall the definition of
quantum affine algebras  and several results on the classification of their
irreducible  finite--dimensional representations.

Let $q$ be an indeterminate, let $\bc(q)$ be the field of rational
functions in $q$ with complex coefficients.  For $r,m\in\bn$, $m\ge r$, define
\begin{equation*}
[m]_q=\frac{q^m -q^{-m}}{q -q^{-1}},\ \ \ \ [m]_q! =[m]_q[m-1]_q\ldots
[2]_q[1]_q,\ \ \ \
\left[\begin{matrix} m\\ r\end{matrix}\right]_q
= \frac{[m]_q!}{[r]_q![m-r]_q!}.
\end{equation*}

 Let $\frak g$ be a complex finite--dimensional simple Lie algebra  of rank $n$,
 set $I=\{1,2,\cdots ,n\}$, let $\{\alpha_i:i\in I\}$ (resp.  $\{\omega_i:i\in I\}$)
  be the set of simple roots
 (resp.  fundamental weights) of $\frak g$ with respect to $\frak h$. As usual, $Q^+$
(resp. $P^+$) denotes the non--negative root (resp. weight)
lattice of $\frak g$.
 Let $A=(a_{ij})_{i,j\in I}$  be the $n\times n$ Cartan matrix of $\frak g$ and
let $\hat A =(a_{ij})$ be  the $(n+1)\times (n+1)$ extended Cartan
matrix associated to $\frak g$. Let $\hat{I} =I\cup\{0\}$. Fix
non--negative integers $d_i$ ($i\in\hat{I}$) such that the matrix
$(d_ia_{ij})$ is symmetric. Set $q_i=q^{d_i}$ and
$[m]_i=[m]_{q_i}$.
\begin{prop}{\label{defnbu}} There is a Hopf algebra $\tilde{\bu}_q$ over
$\bc(q)$ which is generated as an algebra by elements $E_{\alpha_i}$,
$F_{\alpha_i}$, $K_i^{{}\pm 1}$ ($i\in\hat I$), with the following defining
relations:
\begin{align*}
  K_iK_i^{-1}=K_i^{-1}K_i&=1,\ \ \ \ K_iK_j=K_jK_i,\\
  K_iE_{\alpha_j} K_i^{-1}&=q_i^{ a_{ij}}E_{\alpha_j},\\
K_iF_{\alpha_j} K_i^{-1}&=q_i^{-a_{ij}}F_{\alpha_j},\\
  [E_{\alpha_i}, F_{\alpha_j}
]&=\delta_{ij}\frac{K_i-K_i^{-1}}{q_i-q_i^{-1}},\\
  \sum_{r=0}^{1-a_{ij}}(-1)^r\left[\begin{matrix} 1-a_{ij}\\
  r\end{matrix}\right]_i
&(E_{\alpha_i})^rE_{\alpha_j}(E_{\alpha_i})^{1-a_{ij}-r}=0\
  \ \ \ \ \text{if $i\ne j$},\\
\sum_{r=0}^{1-a_{ij}}(-1)^r\left[\begin{matrix} 1-a_{ij}\\
  r\end{matrix}\right]_i
&(F_{\alpha_i})^rF_{\alpha_j}(F_{\alpha_i})^{1-a_{ij}-r}=0\
  \ \ \ \ \text{if $i\ne j$}.
\end{align*}
The comultiplication of $\tilde{\bu}_q$ is given on generators by
$$\Delta(E_{\alpha_i})=E_{\alpha_i}\ot 1+K_i\ot E_{\alpha_i},\ \
\Delta(F_{\alpha_i})=F_{\alpha_i}\ot K_i^{-1} + 1\ot F_{\alpha_i},\ \
\Delta(K_i)=K_i\ot K_i,$$
for $i\in\hat I$.\hfill\qedsymbol
\end{prop}

Set $K_{\theta} =\prod_{i=1}^n K_i^{r_i/d_i}$, where $\theta=\sum
r_i\alpha_i$ is the highest root in $R^+$.
Let $\bu_q$ be the  quotient of $\tilde{\bu}_q$ by the ideal generated by the
central element  $K_0K_{\theta}^{-1}$; we call this the quantum loop algebra
of $\frak g$.

It follows from
\cite{Dr}, \cite{B}, \cite{J} that $\bu_q$ is isomorphic to the
algebra with generators $x_{i,r}^{{}\pm{}}$ ($i\in I$, $r\in\bz$), $K_i^{{}\pm
1}$
($i\in I$), $h_{i,r}$ ($i\in I$, $r\in \bz\backslash\{0\}$) and the following
defining relations:
\begin{align*}
   K_iK_i^{-1} = K_i^{-1}K_i& =1, \ \
 K_iK_j =K_jK_i,\\  K_ih_{j,r}& =h_{j,r}K_i,\\
 K_ix_{j,r}^\pm K_i^{-1} &= q_i^{{}\pm
    a_{ij}}x_{j,r}^{{}\pm{}},\ \ \\
  [h_{i,r},h_{j,s}]=0,\; \; & [h_{i,r} , x_{j,s}^{{}\pm{}}] =
  \pm\frac1r[ra_{ij}]_ix_{j,r+s}^{{}\pm{}},\\
 x_{i,r+1}^{{}\pm{}}x_{j,s}^{{}\pm{}} -q_i^{{}\pm
    a_{ij}}x_{j,s}^{{}\pm{}}x_{i,r+1}^{{}\pm{}} &=q_i^{{}\pm
    a_{ij}}x_{i,r}^{{}\pm{}}x_{j,s+1}^{{}\pm{}}
  -x_{j,s+1}^{{}\pm{}}x_{i,r}^{{}\pm{}},\\ [x_{i,r}^+ ,
  x_{j,s}^-]=\delta_{i,j} & \frac{ \psi_{i,r+s}^+ -
    \psi_{i,r+s}^-}{q_i - q_i^{-1}},\\
\sum_{\pi\in\Sigma_m}\sum_{k=0}^m(-1)^k\left[\begin{matrix}m\\k\end{matrix}
\right]_i
  x_{i, r_{\pi(1)}}^{{}\pm{}}\ldots x_{i,r_{\pi(k)}}^{{}\pm{}} &
  x_{j,s}^{{}\pm{}} x_{i, r_{\pi(k+1)}}^{{}\pm{}}\ldots
  x_{i,r_{\pi(m)}}^{{}\pm{}} =0,\ \ \text{if $i\ne j$},
\end{align*}
for all sequences of integers $r_1,\ldots, r_m$, where $m =1-a_{ij}$, $\Sigma_m$
is the symmetric group on $m$ letters, and the $\psi_{i,r}^{{}\pm{}}$ are
determined by equating powers of $u$ in the formal power series
$$\sum_{r=0}^{\infty}\psi_{i,\pm r}^{{}\pm{}}u^{{}\pm r} = K_i^{{}\pm 1}
{\text{exp}}\left(\pm(q_i-q_i^{-1})\sum_{s=1}^{\infty}h_{i,\pm s} u^{{}\pm
s}\right).$$
\vskip 12pt

For $i\in I$, the preceding  isomorphism maps $E_{\alpha_i}$ to  $x_{i,0}^+$ and
$F_{\alpha_i}$ to $x_{i,0}^-$.  The subalgebra generated by $E_{\alpha_i}$,
$F_{\alpha_i}$, $K_i^{\pm 1}$ ($i\in I$) is the quantized enveloping algebra
$\bu_q^{fin}$ associated to $\frak g$. Let $\bu_q(<)$ be the subalgebra
generated by the elements $x_{i,k}^-$ ($i\in I$, $k\in\bz\}$).  For $i\in I$,
let $\bu_i$ be the subalgebra of $\bu_q$ generated by the elements
$\{x_{i,k}^\pm: k\in\bz\}$, the subalgebra $\bu_i^{fin}$ is defined in the same
way. Notice that $\bu_i$ is isomorphic  to  the quantum affine algebra
$\bu_{q_i}(\widehat{sl_2})$.  Let $\Delta_i$ be the comultiplication of
$\bu_{q_i}(\widehat{sl_2})$.

An explicit formula for the comultiplication on the Drinfeld generators is not
known. Also, the subalgebra $\bu_i$ is not a Hopf subalgebra of $\bu_q$.
However, partial information which is sufficient for our needs is  given in the
next proposition.
Let \begin{equation*} X^\pm =\sum_{i\in I, k\in\bz}  \bc(q)x_{i,k}^\pm,\ \
X^\pm(i) =\sum_{j\in I\backslash\{i\}, k\in\bz}\bc(q)
x_{j,k}^\pm.\end{equation*}

\begin{prop}\label{comultip} The restriction of  $\Delta$ to $\bu_i$ satisfies,
\begin{equation*} \Delta(x) =\Delta_i(x) \mod\left(\bu_q\otimes(\bu_q\backslash
\bu_i)\right).\end{equation*}
More precisely:
\begin{enumerate}
\item[(i)] Modulo $\bu_qX^-\otimes \bu_q(X^+)^2 + \bu_q X^-\otimes \bu_qX^+(i)$,
we have
\begin{align*} \Delta (x_{i,k}^+)& =x_{i,k}^+\otimes 1+ K_i\otimes x_{i,k}^+
+\sum_{j=1}^k \psi^+_{i,j}\otimes x_{i,k-j}^+ \ \ (k\ge 0),\\
\Delta (x_{i,-k}^+)& =K_i^{-1}\otimes x_{i,-k}^+ + x_{i,-k}^+\otimes 1
+\sum_{j=1}^{k-1} \psi^-_{i, -j}\otimes  x_{i,-k+j}^+ \ \ (k> 0).\end{align*}
\item[(ii)]  Modulo $\bu_q(X^-)^2\otimes\bu_qX^+ + \bu_qX^-\otimes \bu_qX^+(i)$,
we have
\begin{align*} \Delta (x_{i,k}^-)& =x_{i,k}^-\otimes K_i+ 1\otimes x_{i,k}^-
+\sum_{j=1}^{k-1} x_{i,k-j}^-\otimes \psi^+_{i, j} \ \ (k>  0),\\
\Delta (x_{i,-k}^-)& = x_{i,-k}^-\otimes K_i^{-1}  + 1\otimes x_{i,-k}^-
+\sum_{j=1}^k  x_{i,-k+j}^-\otimes \psi^-_{i, -j} \ \ (k\ge  0).\end{align*}
\item[(iii)] Modulo $\bu_qX^-\otimes \bu_qX^+$, we have
\begin{equation*} \Delta(h_{i,k}) =h_{i,k}\otimes 1+1\otimes h_{i,k} \ \
(k\in\bz).\end{equation*}
\end{enumerate}
\end{prop}
\begin{pf} Part (iii) was proved in \cite{Da}. The rest of the proposition was
proved in \cite{CPminaff}.
\end{pf}

We conclude this section with some results on the classification of irreducible
finite--dimensional representations of quantum affine algebras.
Let \begin{equation*} \cal{A} =\{f\in\bc(q)[[u]]: f(0)=0\}.\end{equation*}

For any $\bu_q$-module $V$ and any $\mu=\sum_i\mu_i\omega_i\in P$, set
\begin{equation*} V_\mu=\{ v\in V: K_i.v =q_i^{\mu_i}v ,\ \
\forall \ i\in I\}.\end{equation*}
We say that $V$ is a module of type 1 if
\begin{equation*} V=\bigoplus_{\mu\in P}V_\mu.\end{equation*}
From now on, we shall only be working with $\bu_q$-modules of type 1.
For $i\in I$,  set
\begin{equation*}
h^\pm_i(u)=\sum_{k=1}^\infty \frac{q^{\pm k}h_{i,\pm
k}}{[k]_i}u^k.\end{equation*}
 \begin{defn} We say that a $\bu_q$--module $V$ is (pseudo) highest weight, with
highest weight $(\lambda, \boh^\pm)$, where $\lambda=\sum_{i\in
I}\lambda_i\omega_i$,
$\boh^\pm=(\boh_1^\pm(u),\cdots,\boh_n^\pm(u))\in\cal{A}^n$,  if there exists
$0\ne v\in V_\lambda$ such that $V=\bu_q.v$ and
\begin{equation*} x_{i,k}^+. v= 0,\ \ K_i.v =q_i^{\lambda_i}v,\ \
h^\pm_i(u).v=\boh_i^\pm(u)v,\end{equation*}
for all $i\in I$, $k\in\bz$.
\hfill\qedsymbol\end{defn}
If $V$ is any highest weight module, then in fact $V=\bu_q(<).v$ and so
\begin{equation*} V_\mu\ne 0\implies \mu=\lambda-\eta\ \ (\eta\in Q^+).
\end{equation*}
Clearly any highest weight module has  a unique irreducible quotient
$V(\lambda,\boh^\pm)$.

The following was proved in \cite{CP}.
\begin{thm}\label{classify} Assume that the pair $(\lambda, \boh^\pm)\in P\times \cal{A}^n$ satisfies, the following: $\lambda=\sum_{i\in I}\lambda_i\omega_i\in
P^+$, and there exist elements $a_{i,r}\in\bc(q)$ ($1\le r\le \lambda_i, i\in
I$) such that
\begin{equation*}
\boh_i^\pm(u) =-\sum_{r=1}^{\lambda_i}\text{ln}(1-a_{i,r}^{\pm 1}u).
\end{equation*}
Then, $V(\lambda,\boh^\pm)$ is the unique (up to isomorphism)  irreducible
finite--dimensional $\bu_q$--module with highest weight $(\lambda,\boh^\pm)$.
\hfill\qedsymbol\end{thm}
\begin{rem} This statement is actually a reformulation of the statement in
\cite{CP}.  Setting $\pi_i(u)
=\prod_{r=1}^{\lambda_i}(1-a_{i,r}u)$ and calculating the
eigenvalues of the $\psi_{i,k}$ gives the result stated in
\cite{CP}, see also \cite{CPweyl}.
\end{rem}
From now on, we shall only be concerned with the modules
$V(\lambda, \boh^\pm) $ satisfying  the conditions of Theorem
\ref{classify}. In view of the preceding remark, it is clear that
the isomorphism classes of such modules are indexed by an
$n$--tuple of polynomials $\bpi=(\pi_1,\cdots, \pi_n)$, which have
constant term 1, and which are split over $\bc(q)$. We  denote the
corresponding module by $V(\bpi)$ and the highest weight vector by
$v_\bpi$, where $\lambda=\sum_{i=1}^n\text{deg}\pi_i$. For all
$i\in I$, $k\in\bz $, we  have
\begin{equation}\label{qrel1}x_{i,k}^+.v_{\bpi} =0,\ \ K_i.v_{\bpi}
=q_i^{{\text{deg}}\pi_i}v_{\bpi},\end{equation}
and
\begin{equation}\label{qrel2} \frac{q^{\pm k}h_{i,\pm k}}{[k]_i}. v_{\bpi} = \boh_{i, \pm
k}v_{\bpi},\ \
(x_{i,k}^-)^{{\text{deg}}\pi_i+1}.v_{\bpi} =0,\end{equation}
where the $\boh_{i,\pm k}$ are determined from the functional equation
\begin{equation*}{\text {exp}}\left(-\sum_{k>0}{\boh_{i,\pm
k}}u^k\right) = \pi_i^\pm(u),\end{equation*} with
$\pi_i^+(u)=\pi_i(u)$ and $\pi_i^-(u) =
u^{\text{deg}\pi_i}\pi_i(u^{-1})/
\left.\left(u^{\text{deg}\pi_i}\pi_i(u^{-1}\right)\right|_{u=0}$.

For any $\bu_q$--module $V$, let $V^*$ denote its left dual. Let
$- : I\to I$ be the unique  diagram automorphism  such that the
irreducible $\frak g$ module  $V(\omega_i)\cong
V(\omega_{\overline i})$.    There exists an integer $c\in\bz$
depending only on $\frak g$ such that \begin{equation}
V(\bpi)^*\cong V(\bpi^*), \ \ \bpi^*
=(\pi_{\overline{1}}(q^cu),\cdots
,\pi_{\overline{n}}(q^cu)).\end{equation} Analogous statements
hold for right duals \cite{CPminaff}. Recall also, that if a
module and its dual are highest weight then they must be
irreducible.

Finally, let $\omega:\bu_q\to\bu_q$ be the algebra automorphism
and coalgebra anti-automorphism obtained by extending the
assignment $\omega(x_{i,k}^\pm)=-x_{i,-k}^\mp$. If $V$ is any
$\bu_q$-module, let $V^\omega$ be the pull--back of $V$ through
$\omega$. Then, $(V\otimes V')^\omega\cong (V')^\omega\otimes
V^\omega$ and \begin{equation} V(\bpi)^\omega=
V(\bpi^\omega)\end{equation} where
\begin{equation*} \bpi^\omega = (\pi^-_{\overline{1}}(q_1^2\kappa u),\cdots
,\pi^-_{\overline{n}}(q_n^2\kappa u)),\end{equation*} for a  fixed
 $\kappa $ depending only on $\frak g$.

Throughout this paper, we shall only work with polynomials in
$\bc(q)[u]$  which are split and have constant term 1. For any
$0\ne a\in\bc(q)$,  $m\in\bz^+$, set
\begin{equation}\pi_{m,a}(u)=\prod_{r=1}^m (1-aq^{m-2r+1}u).
\end{equation}
   It is a simple combinatorial fact [CP1] that any
such polynomial can be written uniquely as a product
\begin{equation*}\pi(u)=\prod_{j=1}^s\pi_{m_j, a_j},\end{equation*}
where $m_j\in\bz_+ $, $a_j\in\bc(q)$ and
\begin{equation*}j< \ell \implies \frac{a_{j}}{a_\ell }\ne  q^{\pm(m_j+m_{\ell
}-2p)},\ \ 0\le p<\text{min} (m_j,m_\ell).\end{equation*} If $\pi$
and $\pi'$ are two such  polynomials, then  we say that $\pi$ is
in general position with respect to $\pi'$  if for all $1\le j\le
s$, $1\le k\le s'$, we have
\begin{equation*} \frac{a_j}{a'_k}\ne  q^{-(m_j+m_k'-2p)},\end{equation*}
for any  $0\le p < m_j$.  This is equivalent to saying that
 for
all roots  $a$ of  $\pi$ we have,
\begin{equation*}\frac{a}{a_k'}\ne q^{-1-m_k'},\ \forall\  \ 1\le k'\le s'.\end{equation*}

We conclude this section with some results  in the case when $\frak g=sl_2$.

\begin{thm} \label{sl2} \hfill
\begin{enumerate}
\item[(i)]  The irreducible module $V(\pi_{m,a})$ with highest weight $\pi_{m,a}$
is of dimension $m+1$ and is irreducible as a
$\bu_q^{fin}$--module.
\item[(ii)] Assume that $a_1,\cdots , a_\ell\in\bc(q) $ are such that if $r<s$
then $a_r/a_s\ne q^{-2}$ (i.e. the polynomial $(1-a_ru)$ is in
general position with respect to $(1-a_su)$ for all $1\le r\le
s\le \ell$) . The module $W(\pi) = V(\pi_{1,a_1})\otimes\cdots
\otimes V(\pi_{1,a_\ell})$ is the universal finite--dimensional
highest weight module with highest weight
\begin{equation*} \pi(u)=\prod_{r=1}^\ell (1-a_ru),\end{equation*}i.e. any other
finite--dimensional highest weight module with highest weight $\pi$ is a
quotient of $W(\pi)$.
\item[(iii)] For $1\le r\le \ell$, let $a_r\in\bc(q)$ and $m_r\in\bz_+$ be such
that  $\pi_{m_r,a_r}$ is in general position with respect to
$\pi_{m_s, a_s}$ for all $1\le r\le s\le \ell$. The module
$V(\pi_{m_1, a_1}\otimes\cdots\otimes V(\pi_{m_\ell, a_\ell})$ is
a  highest weight module with highest weight $\pi_{m_1,a_1}\cdots
\pi_{m_\ell, a_\ell}$ and   highest weight vector
$v_{\pi_{m_1,a_1}}\otimes \cdots \otimes v_{\pi_{m_\ell,a_\ell}}$.
In particular, the  module $W(\pi_{m_1, a_1})\otimes
V(\pi_{m_2,a_2})\otimes\cdots\otimes V(\pi_{m_\ell, a_\ell})$ is
highest weight.
\item[(iv)] The module $V(\pi)\otimes V(\pi')$ is irreducible iff
\begin{equation*}\frac{a_r}{a_s}\ne q^{\pm (m_r+m_s-2p)}\ \
\forall 0\le p <{\text{min}} (m_r,m_s).\end{equation*}
\end{enumerate}
\end{thm}
\begin{pf} Parts (i)  and (iv)  were   proved  in \cite{CPqa}.  Part (ii) was proved in
\cite{CPweyl}, (see also \cite{VV}, \cite{Nak2} in the case when
$\pi(u)\in\bc[q,q^{-1},u]$. In the general case, choose
$v\in\bc[q,q^{-1}]$ so that $\tilde{\pi}(u)= \pi(uv)$ has all its
roots in $\bc[q,q^{-1}]$. Let $\tau_v:\bu_q\to \bu_q$ be the
algebra and coalgebra automorphism defined by sending $x_k^\pm\to
v^kx_k^\pm$.  The pull back of $V(\pi)$   through $\tau_v$ is
$V(\tilde{\pi})$  and hence $W(\pi(u))\cong W(\pi(uv)$. This
proves (ii). Part (iii) is proved as in  Lemma 4.9 in \cite{CPqa}.
In fact, the proof given there establishes the stronger result
stated here.

\end{pf}

\section{Braid group action} Let $W$ be the Weyl group of $\frak g$ and let
$\cal{B}$ be the corresponding braid group. Thus, $\cal{B}$ is the group
generated by elements $T_i$ ($i\in I$) with defining relations:
\begin{align*} T_iT_j &=T_jT_i,\ \ \text{if}\ \ a_{ij} =0,\\
T_iT_jT_i& =T_jT_iT_j\ \ \text{if}\ \ a_{ij}a_{ji} =1,\\
(T_iT_j)^2&= (T_jT_i)^2,\ \ \text{if}\ \ a_{ij}a_{ji}=2,\\
(T_iT_j)^3&= (T_jT_i)^3,\ \ \text{if}\  \ a_{ij}a_{ji} =3,\end{align*}
where $i,j\in\{1,2,\cdots ,n\}$ and $A=(a_{ij})$ is   the Cartan matrix of
$\frak g$.

A straightforward calculation gives the following proposition.
\begin{prop} For all $r\ge 1$, the formulas
\begin{equation*} T_ie_j =e_j-q_i^r\frac{[ra_{ji}]_j}{[r]_j}e_i,\end{equation*}
define a representation  $\eta_r:\cal{B}\to  \text{end}(V_r)$, where $V_r\cong
\bc(q)^n$ and $\{e_1,\cdots ,e_n\}$ is the standard basis of $V_r$. Further,
identifying
\begin{equation*} \cal{A}^n\cong\prod_{r=1}^\infty V_r,\end{equation*}
we get  a representation of $\cal{B}$ on $\cal{A}^n$ given by
\begin{align*}\label{braid}
(T_i\boh)_j & =\boh_j(u),\ \  {\text{if}}\ a_{ji}=0,\\
(T_i\boh)_j & =\boh_j(u) +\boh_i(qu),\ \  {\text{if}}\ a_{ji}=-1,\\
(T_i\boh)_j & =\boh_j(u) +\boh_i(q^3u) +\boh_i(qu),\ \  {\text{if}}\
a_{ji}=-2,\\
(T_i\boh)_j & =\boh_j(u) +\boh_i(q^5u) +\boh_i(q^3u)+\boh_i(qu),\ \
{\text{if}}\ a_{ji}=-3,\\
 (T_i\boh)_i& =-\boh_i(q_i^2u),\end{align*}
for all $i,j\in I$, $\boh\in\cal{A}^n$.
\hfill\qedsymbol\end{prop}

Let $s_i$, $i\in I$ be a set of simple reflections in $W$. For any
$w\in W$, let $\ell(w)$ be the length of a reduced expression for
$w$.  If $w=s_{i_1}s_{i_2}\cdots s_{i_k}$ is a reduced expression
for $w$, set $I_w=\{i_1,i_2,\cdots i_k\}\subset I$ and
 let $T_w=T_{i_1}\cdots T_{i_k}$.
 It is well--known that $T_w$ and
$I_w$ are  independent of the choice of the reduced expression.
Given $\boh\in\cal{A}^n$ and $w\in W$, we have
\begin{equation*} T_w\boh =T_{i_1}T_{i_2}\cdots T_{i_k}\boh =\left((T_w\boh)_1,
\cdots ,(T_w\boh)_n\right).\end{equation*}
 We can
now prove:
\begin{prop}\label{twh} Suppose that $w\in W$ and $i\in I$ is such that
$\ell(s_iw)=\ell(w)+1$. There exists an integer $M\equiv
M(i,w,\boh)\ge 0$ and non--negative integers $p_{r,j}$ ($j\in I,
1\le r\le M$) such that
 \begin{equation} (T_w\boh)_i=  \sum_{j\in
I_{w}\cup\{i\}}\sum_{r=1}^{M}p_{r,j}\boh_j(q^ru).\end{equation}
Further, if $i\notin I_w$, then \begin{equation} (T_w\boh)_i=
  \boh_i(u)+ \sum_{j\in I_w}
\sum_{r=1}^{M}p_{r,j}\boh_{i_j}(q^ru).\end{equation}
\end{prop}
\begin{pf} Proceed by induction on $\ell(w)$; the induction clearly starts  at
$\ell(w)=0$.
Assume that the result is true for $\ell(w)<k$. If $\ell(w)=k$, write $w=s_jw'$
with $\ell(w')=k-1$. Notice that  $j\ne i$ since $\ell(s_iw)=\ell(w)+1$.
We get
\begin{equation*}(T_w\boh)_i= (T_jT_{w'}\boh)_i
=(T_{w'}\boh)_i(u)+\sum_{s=0}^{|a_{ij}|-1}(T_{w'}\boh)_j(q^{2|a_{ij}|-2s-1}u).
\end{equation*}
If $\ell(s_iw')=\ell(w')+1$, in particular this happens  if $i\not
I_w$,  the proposition follow by induction. If
$\ell(s_iw')=\ell(w')-1$, we have $w=s_js_iw''$. Suppose that
$a_{ij}a_{ji}=-1$. Then, $\ell(s_jw'')=\ell(w'')+1$ and we get
\begin{equation*} (T_jT_iT_{w''}\boh)_i= (T_iT_{w''}\boh)_i+
(T_iT_{w''}\boh)_j(qu) = (T_{w''}\boh)_j(q_iu).\end{equation*}
The result again follows by induction.
The  cases when $a_{ij}a_{ji}=2, 3$ are  proved similarly.  We omit the details.
 \end{pf}

 \section{The main theorem}  Our goal in this section is to obtain a sufficient
condition for a tensor product of two highest weight representations to be
highest weight.

Let $V$ be any highest weight finite--dimensional $\bu_q$--module with highest
weight $\bpi$ (or $(\lambda,\boh^\pm)$ as in Theorem \ref{classify}). For all
$w\in W$, we have
\begin{equation*} \text{dim}\ V_{w\lambda} = 1.\end{equation*}
If  $s_{i_1}\cdots s_{i_k}$ is a reduced expression for $w$, and $\lambda=\sum_i
\lambda_i\omega_i$,  set $m_k =\lambda_k$ and define non--negative integers
$m_j$ (depending on $w$), for $1\le j\le k$, by
\begin{equation*} s_{i_{j+1}}s_{i_{j+2}}\cdots s_{i_k}\lambda
=m_j\omega_j+\sum_{i\ne j}m'_i\omega_i.\end{equation*}
Let $v_\lambda$ be the highest weight vector in $V$. For $w\in W$, set
\begin{equation*}v_{w\lambda} = (x_{i_1,0}^-)^{m_1}\cdots
(x_{i_k,0}^-)^{m_k}.v_\lambda.\end{equation*}
 If $i\in I$ is  such that $\ell(s_iw)=\ell(w)+1$, then
\begin{equation}\label{extr1} x_{i,k}^+. v_{w\lambda} =0,\ \ \forall k\in\bz.
\end{equation} To see this, observe that $w\lambda+\alpha_i$ is not a weight of
$V$, since $w^{-1}\alpha_i\in R^+$ if $\ell(s_iw)=\ell(w)+1$. It is now easy to
see that  $v_{w\lambda}\ne 0$,  $V_{w\lambda} =\bc(q)v_{w\lambda}$ and
\begin{equation*} V=\bu_q.v_{w\lambda}.\end{equation*}
 Since $h_{i,k}V_{w\lambda}\subset V_{w\lambda}$ for all $i\in I$, $k\in\bz$ it
follows that
\begin{equation*} \frac{h_{i,k}}{[k]}.v_{w\lambda} =\boh_{i,k}^wv_{w\lambda}\ \
\forall \ i\in I,\  0\ne k\in \bz,\end{equation*}
where $\boh_{i,k}^w\in\bc(q)$.
Set
\begin{equation*}
\boh^w_i(u) = \sum_{k=1}^\infty \boh^w_{i,k}u^k, \ \ \boh^w= (\boh_1^w(u),
\cdots \boh_n^w(u)).
\end{equation*}
Recall that $\boh^1=-(\text{ln}\pi_1(u),\cdots ,\text{ln}\pi_n(u))$.
\begin{prop} \label{extr} If $w\in W$, then
\begin{equation*} \boh^w = T_w\boh^1.
\end{equation*}
\end{prop}
\begin{pf}
We proceed by induction on $\ell(w)$. If $\ell(w)=0$ then $w=id$ and  the result
follows by definition. Suppose that $\ell(w)=1$, say $w=s_j$. Writing
$\lambda=\sum_j \lambda_j\omega_j$, we have $v_{s_j\lambda}
=(x_{j,0}^-)^{\lambda_j}.v_\lambda$. We first show that
\begin{equation}\label{lw1} h_{j,k}(u). v_{s_j\lambda} = -
\boh_j(q_j^2u)v_{s_j\lambda}  = (T_j\boh)_j(u)v_{s_j\lambda} .\end{equation}
The subspace spanned by the elements $\{(x_{j,l}^-)^r.v_\lambda:0\le r\le
\lambda_j\}$ is a highest weight module for $\bu_j$, hence
it is enough to prove \eqref{lw1}   for highest weight representations of
quantum affine $sl_2$.  In fact it is enough to prove it for the module $W(\pi)$
of Theorem \ref{sl2}. Using Proposition \ref{comultip}, we see  that the
eigenvalue of $h_{i,k}$ on the tensor product of the lowest (and the highest)
weight vectors is
just the sum of the eigenvalues values in  each representation. This reduces us
to the case of  the two-dimensional representation, which is trivial.

Next consider the case $\ell(w)=s_i$, with  $i\ne j$. Recall that \begin{equation*}
[h_{i,r}, x_{j,0}^-] = -\frac{[ra_{ij}]_i}{r} x_{j,r}^- ,\ \ \ \ [h_{j,r},
x_{j,0}^-] =-\frac{[2r]_j}{r}x_{j,r}^-.\end{equation*}
Hence,
\begin{align*} h_{i,r}(x_{j,0}^-)^{\lambda_j} &=
(x_{j,0}^-)^{\lambda_j}h_{i,r}+[ h_{i,r}, (x_{j,0}^-)^{\lambda_j}]\\
&= (x_{j,0}^-)^{\lambda_j}h_{i,r} +\frac{[ra_{ij}]_i}{[2r]_j}[h_{j,r},
(x_{j,0}^-)^{\lambda_j}].\end{align*}
This gives
\begin{align*}\frac{h_{i,r}}{[r]_i}. (x_{j,0}^-)^{\lambda_j}.v_\lambda &=
\boh_{i,r}(x_{j,0}^-)^{\lambda_j}.v_\lambda
+\frac{[ra_{ij}]_i}{(q_j^r+q_j^{-r})[r]_i}\left[\frac{h_{j,r}}{[r]_j},
(x_{i,0}^-)^{\lambda(h_i)}\right].v_\lambda,\\ &=
\boh_{i,r}(x_{j,0}^-)^{\lambda_j}.v_\lambda
+\frac{[ra_{ij}]_{i}}{(q_j^r+q_j^{-r})[r]_i}(-\boh_{j,r}(x_{j,0}^-)^{\lambda_j}.
v_\lambda +\frac {h_{j,r}}{[r]_j}.\left
(x_{j,0}^-)^{\lambda_j}\right).v_\lambda,\\&= \boh_{i,r}
(x_{j,0}^-)^{\lambda_j}.v_\lambda
-\frac{[ra_{ij}]_{i}}{(q_j^r+q_j^{-r})[r]_i}(\boh_{j,r} +
q_j^{2r}\boh_{j,r}). (x_{j,0}^-)^{\lambda_j}).v_\lambda,\\
&=(\boh_{i,r}-q_j^r\frac{[ra_{ij}]_i}{[r]_i}\boh_{j,r})(x_{j,0}^-)^{\lambda_j}.v
_\lambda),\end{align*} where the third equality follows from
\eqref{lw1}.

This proves the result when $\ell(w)=1$. Proceeding by induction on $\ell(w)$,
write $w=s_jw'$ with $\ell(w')=\ell(w)-1$. Since
$v_{w\lambda}=(x_{j,0}^-)^{m_j}v_{w'\lambda}$ for some $m_j\ge 0$, the inductive
step is proved exactly  as in the  case $\ell(w)=1$, with $v_\lambda$ being
replaced by $v_{w'\lambda}$.
 This completes the proof of the proposition.
\end{pf}

Let $w_0\in W$ be the longest element of the Weyl group of $\frak
g$.
\begin{lem}\label{hwsub}
 Let $V$, $V'$ be finite--dimensional highest weight representations with
highest weights $\bpi$ and $\bpi'$  and highest weight vectors $v_{\lambda}$ and
$v_{\lambda'}$ respectively. Assume that $v_{w_0\lambda}\otimes
v_{\lambda'}\in\bu_q(v_\lambda\otimes v_{\lambda'})$. Then, $V\otimes V'$  is
highest weight with highest weight vector $v_\lambda\otimes v_{\lambda'}$ and
highest weight $\bpi\bpi'=(\pi_1\pi_1',\cdots ,\pi_n\pi_n')$.\end{lem}
\begin{pf} It is clear from Proposition \ref{comultip} that the element
$v_\lambda\otimes v_{\lambda'}$ is a highest weight vector with highest weight
$\bpi\bpi'$.  It suffices to prove that
\begin{equation*} V\otimes V'=\bu_q(v_\lambda\otimes
v_{\lambda'}).\end{equation*}
Since $x_{i,k}^-.v_{w_0\lambda} =0$ for all $i\in I$ and $k\in\bz$, it follows
from Proposition \ref{comultip} that
\begin{equation*}\Delta(x_{i,k}^-).(v_{w_0\lambda}\otimes v_{\lambda'})
=v_{w_0\lambda}\otimes x_{i,k}^-. v_{\lambda'}.\end{equation*}
Repeating this argument we see that  that $v_{w_0\lambda}\otimes
V'\subset\bu_q(v_\lambda\otimes v_{\lambda'})$. Now applying  the generators
$E_{\alpha_i}$, $F_{\alpha_i}$ ($i\in\hat{I}$) repeatedly, we see that
 $V\otimes V'\subset \bu_q(v_\lambda\otimes v_{\lambda'})$. This proves the
lemma.
\end{pf}

\begin{lem}\label{slextr} Let $w\in W$ and assume that  $i\in I$ is such that
$\ell(s_iw)=\ell(w)+1$. Then, $v_{w\lambda}\otimes v_{\lambda'}$ generates a
$\bu_i$--highest weight module with highest weight $(T_w\boh)_i\boh_i'$.
\end{lem}

\begin{pf} This is immediate from Proposition \ref{comultip} and Proposition
\ref{extr}. \end{pf}
 We can now prove our main result. Given $\bpi=(\pi_1,\cdots ,\pi_n)$, and $w\in
W$, set
\begin{equation*}T_w\bpi=(\text{exp}-(T_w\ln\pi_1(u))_1, \cdots
,\text{exp}-(T_w\ln\pi_n(u))_n).\end{equation*}

\begin{thm}\label{main} Let $s_{i_1}\cdots s_{i_N}$ be a reduced expression for $w_0$.
The module $V(\bpi_1)\otimes \cdots\otimes V(\bpi_r)$ is highest
weight if for all $1\le j\le N$, $1\le m\le \ell\le r$ the
polynomial $(T_{i_{j+1}}T_{i_{j+2}}\cdots T_{i_N}\bpi_m)_{i_j}$ is
in general position with respect to  $(\bpi_\ell)_{i_j}$.
\end{thm}
\begin{pf} First observe that by Proposition \ref{twh},
 $(T_{i_{j+1}}T_{i_{j+2}}\cdots T_{i_N}\bpi_m)_{i_j}$
is indeed a polynomial. We also claim   that for all $i\in I$, the
polynomial $(\bpi_m)_i$ is in general position with respect to
$(\bpi_\ell)_i$ if $m<\ell$. To see this, choose $1\le j\le N$
maximal so that $i_j=i$, and let $w=s_{i_{j+1}}\cdots s_{i_N}$.
Then $i\notin I_w$ and hence by Proposition \ref{twh}, we have
$$(T_{i_{j+1}}T_{i_{j+2}}\cdots
T_{i_N}\bpi_m)_{i_j}=(\bpi_m)_i(u)\prod_{r=0}^M\prod_{s=j+1}^N((\bpi_m)_s(q^ru))^{p_{r,s}},$$
for some nonnegative integers $M$, $p_{r,s}$. The claim follows
since $(T_{i_{j+1}}T_{i_{j+2}}\cdots T_{i_N}\bpi_m)_{i_j}$ is in
general position with respect to $(\bpi_\ell)_{i_j}$.

If $\frak g=sl_2$, the  theorem was proved in Theorem \ref{sl2}
(ii). For arbitrary $\frak g$, we proceed by induction on $r$. If
$r=1$ there is nothing to prove. Let $r>1$ and let
$V'=V(\bpi_2)\otimes \cdots\otimes V(\bpi_r)$. Then $V'$ is
highest weight module with highest weight vector
$v'=v_{\bpi_2}\otimes\cdots\otimes v_{\bpi_r}$ and highest weight
$\bpi'=\bpi_2\cdots \bpi_r$. Setting $\lambda={\text{deg}} \bpi_1=
((\text{deg}\bpi_1)_1,\cdots ,(\text{deg}\pi_1)_n)$, it is enough
by Lemma \ref{hwsub} to prove that
\begin{equation*} v_{w_0\lambda}\otimes v'\in \bu_q(v_{\bpi_1}\otimes
v').\end{equation*} By  Lemma \ref{slextr}, it suffices to prove
that for all $1\le j\le N$,
\begin{equation*} v_{s_{i_j}s_{i_{j+1}}\cdots s_{i_N}\lambda}\otimes
v'\in\bu_{i_j}(v_{s_{i_{j+1}\cdots s_{i_N}}\lambda}\otimes v').
\end{equation*} In fact it suffices to prove that
\begin{equation*} \bu_{i_j}(v_{s_{i_{j+1}}\cdots s_{i_N}\lambda}\otimes v')=
\bu_{i_j}.v_{s_{i_{j+1}}\cdots s_{i_N}\lambda}\otimes
\bu_{i_j}.v',\end{equation*} as $\bu_{i_j}$--modules.

By Theorem \ref{sl2}(iii), we know that
$\bu_{i_j}.v_{s_{i_{j+1}}\cdots s_{i_N}}$ is a quotient of
$W((T_{i_{j+1}}\cdots T_{i_n}\bpi_1)_{i_j})$. Further we claim
that
\begin{equation*}
\bu_{i_j}.v'=\bu_{i_j}.v_{\bpi_2}\otimes\cdots\otimes\bu_{i_j}.v_{\bpi_r}.
\end{equation*}
To see this, notice that the left hand side is clearly contained
in the right hand side. Since $V'$ is highest weight it follows
that
\begin{equation*}\bu_{i_j}.v_{\bpi_2}\otimes\cdots\otimes\bu_{i_j}.v_{\bpi_r}\subset\bu_q(<).v.
\end{equation*}
Since any element in
$\bu_{i_j}.v_{\bpi_2}\otimes\cdots\otimes\bu_{i_j}.v_{\bpi_r}$ has
weight $\sum_{m=2}^r{\text{deg}}(\bpi_m)-p\alpha_{i_j}$, for some
$p\ge0$, it now follows that
\begin{equation*}\bu_{i_j}.v_{\bpi_2}\otimes\cdots\otimes\bu_{i_j}
.v_{\bpi_r}\subset\bu_{i_j}(<).v' ,
\end{equation*} thus establishing our claim.
Since, $\bu_{i_j}v_{\pi_m)_{i_j}}\simeq V((\bpi_m)_{i_j})$ as
$\bu_q(\hat{sl}_2)$--modules,  we have therefore proved that
$\bu_{i_j}.v_{s_{i_{j+1}\cdots s_{i_N}}\lambda}\otimes
\bu_{i_j}.v'$ is a quotient of the tensor product of
$\bu_{i_j}$--modules  $W((\bpi_1)_{i_j})\otimes
V((\bpi_2)_{i_j})\otimes\cdots\otimes V((\bpi_r)_{i_j})$. Since we
have proved that  the polynomial $(\bpi_m)_{i_j}$ is in general
position with respect to $(\bpi_\ell)_{i_j}$ if $m<\ell$, the
result now follows from Theorem \ref{sl2}(iii).

\end{pf}

\section{Relationship wth Kashiwara's results and conjectures}

Let us consider the special case when $\boh$ has the following form,
\begin{equation*} \boh_j^\pm(u)=0, \ \ j\ne i, \ \ \boh_i^\pm (u)
=-\sum_{r=1}^m\text{ln}(1-aq_i^{m-2r+1}u),\end{equation*}
and denote the corresponding $n$--tuple of power series  by $\boh^i_{m,a}$ and
the $n$--tuple of polynomials by $\bpi^i_{m,a}$.
We shall prove the following result.

\begin{thm}\label{k} Let $k_1,k_2\cdots, k_l\in I$, $a_1,\cdots, a_l\in\bc(q)$,
$m_1,\cdots, m_\ell \in\bz_+$, and assume that
\begin{equation*} r<s\ \ \implies \ \ \frac{a_r}{a_s}\ne
q^{d_{k_r}m_r-d_{k_s}m_s- d_{k_r}-d_{k_s}- p}\ \ \forall\ \ p\ge
0.
\end{equation*}
Then, the tensor product
$V(\bpi^{k_1}_{m_1,a_1})\otimes\cdots\otimes
V(\bpi^{k_\ell}_{m_\ell,a_\ell})$ is a highest weight module.
\end{thm}
Assume the theorem for the moment.

\begin{rem}
In the special case when $m_{j}=1$ for all $j$, it was conjectured in \cite{AK}
that such  a tensor product is cyclic if $a_{j}/a_{l}$ does not have a pole at
$q=0$ if $j<\ell$, and this was proved when $\frak g$ is of type $A_n$ or $C_n$;
subsequently,  a geometric proof of this conjecture  was given in \cite{VV} when
$\frak g$ is simply--laced; a complete proof was given using crystal basis
methods in \cite{K}. \end{rem}

The following corollary to Theorem  \ref{k} was conjectured in
\cite{K}, \cite{HKOTY}.

\begin{cor}
 The  tensor product $V=V(\bpi^{k_1}_{m_1, 1})\otimes\cdots\otimes
V(\bpi^{k_\ell}_{m_\ell,1})$ is an irreducible $\bu_q$--module.
\end{cor}
\begin{pf} First observe that if $d_{k_1}m_1\le d_{k_2}m_2\le \cdots \le
d_{k_l}m_\ell$ then the tensor product is cyclic by Theorem
\ref{k}. We claim that it suffices to prove the corollary in the
case when $\ell=2$.  For then, by rearranging the factors in the
tensor product we can show that both $V$ and its dual are  highest
weight and hence  irreducible.  To see that $V=V(\bpi^{k_1}_{m_1,
1})\otimes V(\bpi^{k_2}_{m_2, 2})$ is cyclic if $d_{k_1}m_1>
d_{k_2}m_2$, we have to consider the case when
$d_{k_1}m_1-d_{k_2}m_2-d_{k_1}-d_{k_2}\ge 0$. It suffices to show
that $V^\omega$ is cyclic, since $\omega$ is an algebra
automorphism. Now, $V^\omega = V(\bpi^{\overline{k}_2}_{m_2,
q_2^2\kappa})\otimes V(\bpi^{\overline{k}_1}_{m_1, q_1^2\kappa})$
for some fixed $\kappa$ depending only on $\frak g$. By Theorem
\ref{k}, this is cyclic, since $2d_{k_2}-2d_{k_1}\ne d_{k_2}m_2-
d_{k_1}m_1- d_{k_1}-d_{k_2}-p$\  for any $p\ge 0$.
 This proves
the result.

 \end{pf}

It remains to prove the theorem, for which we must show that,  if
$w_0=s_{i_1}\cdots s_{i_N}$ and $w=s_{i_{j+1}}\cdots s_{i_N}$,
 then the polynomial
$\left(T_w\bpi^{k_r}_{m_r,a_r}\right)_{i_j} $ is in general
position with respect to $ (\bpi^{k_\ell}_{m_\ell,a_\ell})_{i_j}$
for all $r<\ell$.

Using Propositon \ref{twh}, we see that
\begin{equation*} \left(T_w\bpi^{k_r}_{m_r,a_r}\right)_{i_j} =\prod_{s\ge
0}\pi_{m_r,a_r}(q^su),\end{equation*} where $s$ varies over a
finite subset of $\bz_+$ with multiplicity. This means that  any
root of $(T_w\bpi^{k_r}_{m_r,a_r})_{i_j}$ has the form
$q^{d_{k_r}(m_r-2p+1)+s}a_r$ where $s\ge 0$ and $1\le p\le m_r$,
If $\ell\ne i_j$, there is nothing to prove since $
(\bpi^{k_\ell}_{m_\ell,a_\ell})_{i_j}=1$. If $\ell=i_j$, then
 the assumption on $a_r/a_\ell$ implies  that
\begin{equation*}  \frac{q^{d_{k_r}(-m_r+2p+1)+s}a_r}{a_l}\ne
q^{-d_{k_\ell}(1+m_\ell)}.\end{equation*} This   proves the
theorem.

\section{The cyclicity condition made explicit}
In this section we work with a specific reduced expression for the
longest element $w_0$  of the Weyl group,  and give  the condition
explicitly for the tensor product $V(\bpi)\otimes V(\tilde\bpi)$
to be cyclic. In what follows we assume that the nodes of the
Dynkin diagram of $\frak g$ are numbered so that 1 is the short
root (resp. long root) when $\frak g=B_n$ (resp. $\frak g=C_n$)
and that 1 and 2 are the spin nodes  when $\frak g=D_n$.

\subsection{The classical algebras} The main result is:
\begin{thm}\label{explicit}
\begin{enumerate}
\item[(i)] Assume that $\frak g$ is of type $A_n$. Then
 $V(\bpi)\otimes V(\bpi')$ is cyclic if  for all $1\le j\le
i\le n$, we have that the polynomial   $$
\prod_{r=i-j+1}^i\pi_r(q^{2i-j-r}u)\pi_{i-j+2}(q^{i-2}u)\cdots
\pi_i(q^{i-j}u)$$ is in general position with respect to $\pi_j'$.
\item[(ii)] Assume that $\frak g$ is of type $B_n$. Then  $V(\bpi)\otimes V(\tilde\bpi)$ is
cyclic if for all $i\ge 1$ we have that the polynomial $$
\pi_1(q^{4i-4}u)\prod_{r=2}^i\pi_r(q^{4i-1-2r}u)\pi_r(q^{4i-3-2r}u)
 $$ is in general position with respect to $\pi_1'$  and for all $i\ge j\ge
 2$,the polynomials
\begin{eqnarray*}\pi_1(q^{4i-2j-3}u)\prod_{r=2}^i\pi_r(q^{4i-2j-2r}u)\prod_{r=2}^{j-1}&
\pi_r(q^{4i-2j-6+2r}u),\\ \prod_{r=j}^i\pi_r(q^{4i-6+2j-2r}u)&
\end{eqnarray*} are in general position with respect to $\pi_j'$.
\item[(iii)] Assume that $\frak g$ is of type $C_n$. Then  $V(\bpi)\otimes V(\tilde\bpi)$ is
cyclic if for all $i\ge 1$ we have that the polynomial  $$
\prod_{r=1}^i\pi_r(q^{2i-1-r}u)$$ is in general position with
respect to $\pi_1'$   and for all $i\ge j\ge 2$ the polynomials
\begin{eqnarray*}&\pi_1(q^{2i-j-1}u)\pi_1(q^{2i-j+1}u)
\prod_{r=2}^i\pi_r(q^{2i-j-r}u)\prod_{r=2}^{j-1}\pi_r(q^{2i-j+r}u)
\\ & \prod_{r=j}^i\pi_r(q^{2i-r+j}u) \end{eqnarray*} are  in
general position with respect to $\pi_j'$
\item[(iv)] Assume that $\frak g$ is of type $D_n$. Then  $V(\bpi)\otimes V(\tilde\bpi)$ is
cyclic if for all $i\ge 3$, and $i$ even (resp. $i$ odd), we have
\begin{eqnarray*}
\pi_1(q^{2i-4}u)\prod_{r=3}^i\pi_r(q^{2i-2-r}u),\\ (resp.\
\pi_2(q^{2i-4}u)\prod_{r=3}^i\pi_r(q^{2i-2-r}u))\end{eqnarray*} is
in general position with respect to $\pi_1'$,
\begin{eqnarray*}
\pi_2(q^{2i-4}u)\prod_{r=3}^i\pi_r(q^{2i-2-r}u),\ \ (resp.\
\pi_1(q^{2i-4}u)\prod_{r=3}^i\pi_r(q^{2i-2-r}u)),\end{eqnarray*}
is in general position with respect to $\pi_2'$ and  and for all
$i\ge j\ge 3$, the polynomials
\begin{eqnarray*}
\pi_1(q^{2i-j-2}u)\pi_2(q^{2i-j-2}u)\prod_{r=3}^i\pi_r(q^{2i-j-r}u)
\prod_{r=3}^{j-1}\pi_r(q^{2i-j+r-4}u,&\\
\prod_{r=j}^i\pi_r(q^{2i-r+j-4}u),&
\end{eqnarray*}
are in general position with respect to $\pi_j'$.
\end{enumerate}
\end{thm}

We note the following corollary  of this theorem which gives us a
 condition for a  tensor product  $V(\bpi^{k_1}_{m_1,a_1})\otimes\cdots\otimes
V(\bpi^{k_\ell}_{m_\ell,a_\ell})$ to be cyclic. It suffices in
view of Theorem \ref{main} to give the condition for a pair of
such representations to be cyclic. Recall also that the
representation $V(\bpi^{i_1}_{m_1,a_1})\otimes
V(\bpi^{i_2}_{m_2,a_2})$ is cyclic if and only if the
representation $V(\bpi^{i_2}_{m_2,q_{i_2}^2a_2^{-1}})\otimes
V(\bpi^{i_1}_{m_1, q_{i_1}^{-1}a_1^{-1}})$ is cyclic.

\begin{cor} The tensor product $V(\bpi^{i_1}_{m_1,a_1})\otimes V(\bpi^{i_2}_{m_2,a_2})$ is cyclic if
$a_1^{-1}a_2\notin {\cal{S}}(\bpi^{i_1}_{m_1,a_1},
\bpi^{i_2}_{m_2,a_2})$, where $$ {\cal{S}}(\bpi^{i_1}_{m_1,a_1},
\bpi^{i_2}_{m_2,a_2})=\bigcup_{p=1}^{{\text{min}}(m_1,m_2)}q_{i_2}^{d_{i_1}m_1+d_{i_2}m_2+d_1+d_2
-2p}{\cal{S}}(i_1,i_2)$$ and ${\cal{S}}(i_1,i_2)$ is defined as
follows:
\begin{enumerate}
\item[(i)] $\frak g=A_n$, $${\cal{S}}(i_1,i_2)=\{q^{2k-i-j}:1\le i_1,i_2\le k\le n,\   i_1+i_2
\ge k+1\}.$$
\item[(ii)] $\frak g=B_n$,
\begin{eqnarray*}& {\cal{S}}(1,1)=\{q^{4k-4}: 1\le k\le n\},\\
&{\cal{S}}(1,i_2)=\{q^{4k-2i_2-3}: 2\le i_2\le n\},\ \\ &{\cal
S}(i_1,1)=\{q^{4k-2i_1-3}, \ q^{4k-2i_1-1}: 2\le i_1\le n\}, \\
&{\cal S}(i_1,i_2)=\{ q^{4k-2i_2-2i_1}, \ \ q^{4k-2{\text{
max}}(i_1,i_2)+2\ {\text {min}}(i_1,i_2)-6}: 2\le i_1, i_2\le k\le
n,\}.\end{eqnarray*}

\item[(iii)] $\frak g=C_n$,
\begin{eqnarray*}
&{\cal{S}}(i_1,1)=\{q^{2k-1-i_1}: 1\le i_1 \le k\le n \},\\
&{\cal{S}}(1,i_1)=\{q^{2k-1-i_1},\ q^{2k+1-i_1} : 1\le i_1 \le
k\le n \},\\ &{\cal{S}}(i_1,i_2)=\{ q^{2k-i_1-i_2}, \ \
q^{2k-{\text {max}}(i_1,i_2)+{\text{ min}}(i_1,i_2)}:2\le i_1,
i_2\le k\le n\}.\end{eqnarray*}
\item[(iv)] $\frak g =D_n$,
\begin{eqnarray*}& {\cal{S}}(1,1)=\{q^{2k-4}:1\le k\le n,
k\equiv 0\mod 2\} ={\cal S}(2,2),\\
&{\cal{S}}(1,2)=\{q^{2k-4}:1\le k\le n, k\equiv 1\mod 2\}\\ &
{\cal{S}}(i_1,1)=\{q^{2k-i_1-2}:3\le i_1\le k\le
n\}={\cal{S}}(i_1,2) ,\\ & {\cal{S}}(i_1,i_2)=\{q^{2k-i_1-i_2},\ \
q^{2k-{\text{max}}(i_1,i_2)+{\text {min}}(i_1,i_2)-4}: 3\le i_1,
i_2 \le k\le n\}.\end{eqnarray*}
\end{enumerate}
\end{cor}

To prove Theorem \ref{explicit}, we fix a reduced expression for
the longest element $w_0$. Thus, we take
$w_0=\gamma_n\gamma_{n-1}\cdots\gamma_1$, where
\begin{eqnarray*}\gamma_i= \begin{cases}s_1s_2\cdots s_i,\ \ \frak g=A_n,\\
s_is_{i-1}\cdots s_1s_2\cdots s_i,\ \ \frak g=B_n,C_n,\\
s_is_{i-1}\cdots s_2s_1s_3\cdots s_i,\ \ \frak
g=D_n.\end{cases}\end{eqnarray*} It is convenient also to set $w_i
= \gamma_{n-i}\gamma_{n-i-1}\cdots \gamma_1$ for $0\le i\le n-1$.

Theorem \ref{explicit} is an immediate consequence of the
following proposition, which is easily established by an induction
on $i$. We state the proposition only for algebras of type $A_n$
and $B_n$, the results for the other algebras are entirely
similar.

\begin{prop} Let $\boh\in{\cal A}^n$.
\begin{enumerate}
\item[(i)] $\frak g=A_n$. For all $i\ge j$, we have
 \begin{eqnarray*} & (T_{j+1}\cdots
T_{i}T_{w_{i-1}}\boh)_j =\sum_{r=i-j+1}^i\boh_r(q^{2i-j-r}u),\\
&T_{w_i}(\boh)=(-\boh_i(q^{i+1}u),\cdots, -\boh_1(q^{i+1}u),\
\sum_{r=1}^{i+1}\boh_r(q^{i+1-r}u),\ \boh_{i+2}(u), \cdots ,
\boh_n(u)).\end{eqnarray*}
\item[(ii)]$\frak g=B_n$. For all $i\ge j$,
\begin{eqnarray*}
(T_{j+1}T_{j+2}\cdots T_iT_{w_{i-1}}\boh)_j =
&\boh_1(q^{4i-2j-3}u)+\sum_{r=2}^{i}\boh_r(q^{{4i-2j-2r}}u)
+\sum_{r=2}^{j-1}\boh_r(q^{4i-2j-6+2r}u),\\ (T_2T_3\cdots
T_iT_{w_{i-1}}\boh)_1
=&\boh_1(q^{4i-4}u)+\sum_{r=2}^i\boh_r(q^{4i-1-2r}u)+\boh_r(q^{4i-3-2r}u),\\
(T_{j-1}\cdots T_2T_1\cdots T_iT_{w_{i-1}}\boh)_j =&
\sum_{r=j}^i\boh_r(q^{4i+2j-6-2r}u),\end{eqnarray*}
\begin{equation*}(T_{w_i}\boh)_j=\begin{cases} -\boh_j(q^{4i-2}u) & if \ j\le i,\\
=\boh_1(q^{2i-1}u)+\sum_{r=2}^{i+1}\boh_r(q^{{2i-2r+2}}u)
+\sum_{r=2}^{j-1}\boh_r(q^{4i+2r-4}u) & if \ j=i+1,\\
 =\boh_j(u) & if j>i+2.\end{cases}\end{equation*}
\end{enumerate}

\end{prop}

\subsection{The Exceptional algebras} We content ourselves with writing down
 the set ${\cal S}(i_1,i_2)$,
$i_1\le i_2$ which gives the values of $a_1^{-1}a_2$ for which the
tensor product $V(i_1,a_1)\otimes V(i_2,a_2)$ is not cyclic, as
usual if $i_1>i_2$ the condition is that
$q_1^2q_2^{-2}a_1^{-1}a_2\notin {\cal S}(i_2,i_1)$ . The values
were obtained by using mathematica, the program can be used to
write down the conditions for an arbitrary tensor product of
representations  to be irreducible.

We first consider the algebras $E_n$, $n=6,7,8$, the nodes are
numbered as in \cite{Bo} so that 2 is the special node.  The
reduced expression for the longest element is chosen as in
\cite{CXi}, namely:
\begin{eqnarray*} E_6 \ \
\ s_1s_3s_4s_2s_5s_4s_3s_1s_6s_5s_4s_2s_3s_4s_5s_6u_1,\\ E_7\ \ \
s_7s_6s_5s_4s_2s_3s_4s_5s_6s_7s_1s_3s_4s_5
s_2s_4s_3s_1s_6s_5s_4s_2s_3s_4s_5s_6s_7u_2,\\ E_8\ \ \
s_8s_7s_6s_5s_4s_2s_3s_4s_5s_6s_7s_8s_1s_3s_4s_2s_2
s_5s_4s_3s_1s_6s_5s_7s_6s_4s_3s_2s_5s_4s_5s_2s_3s_4\\ \times
s_6s_5s_7s_6s_1s_3s_4s_2s_5s_4s_3s_1s_8s_7s_6s_5s_4s_2s_3s_4s_5s_6s_7s_8u_3,\end{eqnarray*}
where $u_1$ (resp. $u_2$, $u_3$) are the longest elements of $D_5$
obtained by dropping the node 6 (resp. 7,8) (resp. $E_6$, $E_7$)
chosen previously.

\begin{prop}
\begin{enumerate}
\item[(i)] Assume that $\frak g=E_6$. Then,
\begin{eqnarray*}
&{\cal S}_6(1,1) =\{q^2, q^8\} = q^{-4}{\cal  S}_6(1,6),\\
 &{\cal S}_6(1,2)=\{q^5, q^9\},\\
 & {\cal S}_6(1,3)=\{ q^3,q^7, q^9\},\\
& {\cal S}_6(1,4) =\{q^4,q^6,q^8,q^{10}\},\\ & {\cal S}_6(1,5)=
\{q^5, q^7, q^{11}\},\\
\\
&{\cal S}_6(2,2) =\{q^2, q^6, q^8, q^{12}\},\\ & {\cal
S}_6(2,3)=\{q^4,q^6,q^8,q^{10}\}=  {\cal S}_6(2,5),\\ &{\cal
S}_6(2,4)=\{q^3, q^5, q^7,q^9,q^{11}\},\\ & {\cal S}_6(2,6)=\{q^5,
q^9\},\\ \\ & {\cal S}_6(3,3)=\{q^2, q^4,q^6, q^8,q^{10}\},\\ &
{\cal S}_6(3,4)=\{q^3, q^5, q^7, q^9, q^{11}\},\\ & {\cal
S}_6(3,5)=\{q^4, q^6, q^8, q^{10}, q^{12} \},\\
 &  {\cal S}_6(3,6)=\{q^5, q^7,q^{11}\}\\
   \\
   & {\cal S}_6(4,4)= \{q^2,q^4, q^6, q^8, q^{10}, q^{12} \},\\
   &  {\cal S}_6(4,5)=\{q^3, q^5, q^7, q^9, q^{11}\},\\
& {\cal S}_6(4,6)=\{q^4, q^6, q^8,  q^{10}\},\\
\\
&{\cal S}_6(5,5)=\{q^2, q^4,q^6, q^8,q^{10}\},\\
 & {\cal S}_6(5,6) =\{ q^3,q^7, q^9\},\\ \\ & {\cal S}_6(6,6) =\{q^2,
q^8\}.\end{eqnarray*}
\item[(ii)] If $\frak g=E_7$, then ${\cal S}(i,j)$is the union of
the set ${\cal S}_6(i,j)$ with the sets ${\cal S}_7(i,j)$ defined
below.
\begin{eqnarray*}
&{\cal S}_7(1,1) =\{q^{12}, q^{18}\},\\
 &{\cal S}_7(1,2)=\{q^{11}, q^{13}\},\\
 & {\cal S}_7(1,3)=\{ q^{11},q^{13}, q^{15}\},\\
& {\cal S}_7(1,4) =\{q^{12},q^{14},q^{16}\},\\ & {\cal S}_7(1,5)=
\{q^9, q^{13}, q^{15}\},\\
 & {\cal  S}_7(1,6)=\{q^8, q^{14}\}\\
 & {\cal S}_7(1,7)=\{q^7\}\end{eqnarray*}
\begin{eqnarray*}
&{\cal S}_7(2,2) =\{q^{10}, q^{14}, q^{18}\},\\
 & {\cal
S}_7(2,3)=\{q^{12},q^{16}\} \\ &{\cal S}_7(2,4)=\{q^9,q^{13},
q^{15}, q^{17} \},\\ &{\cal S}_7(2,5)=\{q^8, q^{12}, q^{14},
q^{16}\}\\
 & {\cal S}_7(2,6)=\{q^7,
q^{11}, q^{13}, q^{15}\},\\& {\cal S}_7(2,7)=\{q^6, q^{10},
q^{14}\} \\ \\ & {\cal S}_7(3,3)=\{q^{12}, q^{14},q^{16},
q^{18}\},\\ & {\cal S}_7(3,4)=\{q^{13}, q^{15}, q^{17}\},\\ &
{\cal S}_7(3,5)=\{q^{14}, q^{16} \},\\
 &  {\cal S}_7(3,6)=\{q^9, q^{13},q^{15}\}\\
&{\cal S}_7(3,7)=\{q^6, q^8,q^{12}, q^{14}\}\\
   \\
   & {\cal S}_7(4,4)= \{q^{14}, q^{16}, q^{18} \},\\
   &  {\cal S}_7(4,5)= \{q^{13}, q^{15}, q^{17}\},\\
& {\cal S}_7(4,6)=\{q^{12}, q^{14}, q^{16}\},\\
 &{\cal
S}_7(4,7)=\{q^5, q^7, q^9, q^{11}, q^{13}, q^{15}\}
\\
\\
&{\cal S}_7(5,5)=\{q^{12}, q^{14},q^{16}, q^{18}\},\\
 & {\cal S}_7(5,6) =\{ q^5,q^{11}, q^{13}, q^{15}, q^{17}\},\\
& {\cal S}_7(5,7)=\{q^4, q^8, q^{10}, q^{12}, q^{16}\} \\ \\ &
{\cal S}_7(6,6) =\{q^4, q^{10}, q^{12}, q^{16}, q^{18}\},\\ &{\cal
S}_7(6,7)=\{q^3, q^9, q^{11}, q^{17}\}\\ &{\cal S}_7(7,7)=\{q^2,
q^{10}, q^{18}\} .\end{eqnarray*}
 \item[(iii)]If $\frak g=E_8$, then ${\cal S}(i,j)$is the
union of the set ${\cal S}_6(i,j)\cup {\cal S}_7(i,j)$ with the
sets ${\cal S}_8(i,j)$ defined below.
\begin{eqnarray*}
&{\cal S}_8(1,1) =\{q^{14}, q^{20}, q^{24}, q^{30}\},\\
 &{\cal S}_8(1,2)=\{q^{17}, q^{19}, q^{21}, q^{23}, q^{27}\},\\
 & {\cal S}_8(1,3)=\{ q^{17},q^{19}, q^{21}, q^{23}, q^{25},q^{29}\},\\
& {\cal S}_8(1,4) =\{q^{18},q^{20},q^{22}, q^{24}, q^{26},
q^{28}\},\\ & {\cal S}_7(1,5)= \{q^{11}, q^{17}, q^{19}, q^{21},
q^{23},q^{25}, q^{27}  \},\\
 & {\cal  S}_8(1,6)=\{q^{10}, q^{16}, q^{18}, q^{20}, q^{24}, q^{26} \}\\
 & {\cal S}_8(1,7)=\{q^9, q^{17}, q^{19}, q^{25}\}\\
 & {\cal S}_8(1,8)=\{q^8,q^{14}, q^{18}\}
\\
&{\cal S}_8(2,2) =\{q^{12}, q^{16}, q^{20}, q^{22}, q^{24},
q^{26}, q^{30}\},\\
 & {\cal
S}_8(2,3)=\{q^{14},q^{18}, q^{20}, q^{22}, q^{24}, q^{26}, q^{28}
\}
\\ &{\cal S}_8(2,4)=\{ q^{19}, q^{21},
 q^{23}, q^{25}, q^{27}, q^{29} \},\\ &{\cal
S}_8(2,5)=\{q^{18}, q^{20}, q^{22}, q^{24}, q^{26}, q^{28}\}\\
 & {\cal S}_8(2,6)=\{q^{17},
q^{19}, q^{21}, q^{23}, q^{25}, q^{7}\},\\& {\cal
S}_8(2,7)=\{q^{12}, q^{16}, q^{18}, q^{20}, q^{22},q^{24}, q^{26},
q^{28} \}\\
 & {\cal S}_8(2,8)=\{q^7, q^{11}, q^{17}, q^{21}, q^{25}\} \\
\\ & {\cal S}_8(3,3)=\{q^{20}, q^{22},q^{24}, q^{26}, q^{28}, q^{30}\},\\ &
{\cal S}_8(3,4)=\{q^{19}, q^{21}, q^{23}, q^{25}, q^{27},
q^{29}\},\\ & {\cal S}_8(3,5)=\{q^{18}, q^{20}, q^{22}, q^{24},
q^{26}, q^{28} \},\\
 &  {\cal S}_8(3,6)=\{q^{17}, q^{19},q^{21}, q^{23}, q^{25}, q^{27}\}\\
&{\cal S}_8(3,7)=\{q^{10}, q^{16}, q^{18}, q^{20}, q^{22}, q^{24},
q^{26} \}\\ &{\cal S}_8(3,8) =\{ q^7, q^9, q^{13}, q^{15}, q^{17},
q^{19}, q^{23}, q^{25}\}\\
   \\
   & {\cal S}_8(4,4)= \{q^{20}, q^{22}, q^{24}, q^{26}, q^{28}, q^{30} \},\\
   &  {\cal S}_8(4,5)= \{q^{19}, q^{21}, q^{23}, q^{25}, q^{27}\},\\
& {\cal S}_8(4,6)=\{q^{18}, q^{20}, q^{22}, q^{24}, q^{26},
q^{28}\},\\
 &{\cal
S}_8(4,7)=\{q^{17}, q^{19}, q^{21}, q^{23}, q^{25}, q^{27}\}
\\ & {\cal S}_8(4,8)=\{ q^6, q^8, q^{10}, q^{12}, q^{14}, q^{16},
q^{18},q^{20}, q^{22}, q^{24}, q^{26}\}
\\ \\
&{\cal S}_8(5,5)=\{q^{20}, q^{22},q^{24}, q^{26}, q^{28}, q^{30}
\},\\
 & {\cal S}_8(5,6) =\{q^{19}, q^{21}, q^{23}, q^{25},
  q^{27}, q^{29} \},\\
& {\cal S}_8(5,7)=\{q^6, q^{18}, q^{20}, q^{22}, q^{24}, q^{26},
q^{28} \}
\\ & {\cal S}_8(5,8)=\{q^5, q^9, q^{11}, q^{13}, q^{15}, q^{17}, q^{19}\}\\
 \\
 & {\cal S}_8(6,6) =\{q^6, q^{14}, q^{20}, q^{22},
q^{24}, q^{26}, q^{28}, q^{30}\},\\ &{\cal S}_8(6,7)=\{q^5,
q^{13}, q^{15}, q^{19}, q^{21}, q^{23}, q^{27},q^{29} \}\\ & {\cal
S}_8(6,8)=\{q^4, q^8, q^{12}, q^{14}m q^{18}, q^{20}, q^{22},
q^{28}\}
\\ &{\cal S}_8(7,7)=\{q^{12}, q^{14}, q^{20}, q^{22}, q^{28}, q^{30}\}\\
& {\cal S}_8(7,8)=\{q^{3}, q^{11}, q^{13}, q^{19}. q^{21},
q^{29}\},
\\
\\
 &{\cal S}_8(8,8)=\{q^2, q^{12}, q^{20}, q^{30}\}
.\end{eqnarray*}
\end{enumerate}
\end{prop}

If $\frak g=F_4$, then we work with the following reduced
expression for $w_0$, \begin{equation*}
s_4s_3s_2s_3s_4s_1s_2s_3s_2s_1s_4s_3s_2s_3s_4s_1s_2s_3s_2s_1s_2s_3s_2s_3.\end{equation*}
We assume here that the nodes 1 and 2 correspond to the short
simple roots. Assume that $d_i\ge d_j$, then  $V(i,a)\otimes
V(j,b)$ is cyclic if $ab^{-1}\in{\cal S}(i,j)$, where
\begin{eqnarray*} &{\cal S}(1,1)= \{q^2, q^8, q^{12}, q^{18}\},\\
&{\cal S}(1,2)= \{q^7, q^9,   q^{13}, q^{17}\}, \\& {\cal
S}(1,3)=\{q^6, q^{10}, q^{16}\},\\ & {\cal S}(1,4))=\{ q^8,
q^{14}\}.\\
\\ & {\cal S}(2,2)=\{q^2,  q^6, q^8, q^{10}, q^{12},
q^{14}, q^{16}, q^{18}\},\\ & {\cal S}(2,3)= \{q^5, q^7, q^{9},
q^{11}, q^{13}, q^{15}, q^{17}\},
\\
& {\cal S}(2,4)=\{q^7, q^{11}, q^{13}, q^{15}\},
\\
\\
& {\cal S}(3,3)=\{ q^6, q^8, q^{10}, q^{12}, q^{14}, q^{16},
q^{18}\},\\
 & {\cal S}(3,4)= \{q^6, q^8, q^{10},  q^{14}, q^{16}\},\\
\\
\\
& {\cal S}(4,4)=\{q^4, q^{10},q^{12}, q^{18}\} \end{eqnarray*}

If $\frak g =G_2$, then we choose the following reduced expression
for $w_0$, where 1 is corresponds to the short simple root.
\begin{equation*} s_2s_1s_2s_1s_2s_1.\end{equation*}
We have,
\begin{eqnarray*}& {\cal S}(1,1)= \{ q^2, q^6, q^8, q^{12}\}\\
& {\cal S}(1,2)=\{ q^7, q^{11}\} \\ \\& {\cal S}(2,2)=\{ q^6,
q^{8}, q^{10}, q^{12} \}
\end{eqnarray*}

\end{document}